\documentclass[12pt,twoside]{article}  

\usepackage{amsmath}
\usepackage{amsthm}
\usepackage{amssymb}
\usepackage{amscd}
\usepackage{amsfonts}
\usepackage{dsfont}
\usepackage{a4}
\usepackage[dvipsnames,usenames]{color}
\usepackage{eurosym}
\usepackage{fancyhdr}
\usepackage{graphicx}
\usepackage{graphics}
\usepackage[iso-8859-7]{inputenc}
\usepackage{latexsym}
\usepackage{makeidx}
\usepackage{appendix}
\usepackage[dvipsnames,usenames]{color}
\DeclareGraphicsExtensions{.pdf,.jpg}
\usepackage[all]{xy}
\input xy
\xyoption{all}

\newtheorem{theo}{\small\bf Theorem}[section]
\newtheorem{lem}{\small\bf Lemma}[section]

\newtheorem{rem}{\small\bf Remark}[section]

\newcommand{\law}{\stackrel{\mbox{\footnotesize d}}{=}}
\newcommand{\lawv}{\stackrel{\mbox{\footnotesize{\rm d}}}{\longrightarrow}}
\newcommand{\essinf}{{\rm \hspace*{.2ex}ess\inf}}
\newcommand{\esssup}{{\rm \hspace*{.2ex}ess\sup}}
\newcommand{\Z}{\mathds{Z}}
\newcommand{\R}{\mathds{R}}
\newcommand{\E}{\mathds{E}}
\newcommand{\Var}{\mbox{\rm \hspace*{.2ex}Var\hspace*{.2ex}}}
\newcommand{\MATRIXS}{\bbb{\Sigma}\hspace{-1.25ex}|\hspace{-.6ex}|\hspace{.57ex}}

\renewcommand{\d}{\hspace{.2ex}\textrm{$\rm{d}$}}
\newcommand{\ds}{\displaystyle}

\newcommand{\CircleB}{\mbox{\setlength{\unitlength}{.5mm}
 \begin{picture}(2,2)(0,-1.6)\put(0,0){\circle*{2}}\end{picture}}}
\newcommand{\Circle}{\mbox{\setlength{\unitlength}{.5mm}
 \begin{picture}(2,2)(0,-1.6)\put(0,0){\circle{2}}\end{picture}}}
\newcommand{\SquareB}{\mbox{\setlength{\unitlength}{.5mm}
 \begin{picture}(2,2)(0,-1.6)\multiput(-1,-1)(0,.1){20}{\line(1,0){2}}\end{picture}}}
\newcommand{\Square}{\mbox{\setlength{\unitlength}{.5mm}
 \begin{picture}(2,2)(0,-1.6)\put(-1,-1){\line(1,0){2}}\put(-1,-1){\line(0,1){2}}\put(-1,1){\line(1,0){2}}\put(1,-1){\line(0,1){2}}\end{picture}}}
\newcommand{\TriangleB}{\mbox{\setlength{\unitlength}{.35mm}
 \begin{picture}(3,3)(0,-2.3)\put(-2.3,-.5){\scriptsize$_{\blacktriangle}$}\end{picture}}}
\newcommand{\Triangle}{\mbox{\setlength{\unitlength}{.35mm}
 \begin{picture}(3,3)(0,-2.3)\put(-2.3,-.5){\scriptsize$_{\vartriangle}$}\end{picture}}}
\newcommand{\DiamondB}{\mbox{\setlength{\unitlength}{.35mm}
 \begin{picture}(3,3)(0,-2.3)\put(-2.3,1.2){\tiny$_{_\blacklozenge}$}\end{picture}}}
\renewcommand{\Diamond}{\mbox{\setlength{\unitlength}{.35mm}
 \begin{picture}(3,3)(0,-2.3)\put(-2.3,1.2){\tiny$_{_\lozenge}$}\end{picture}}}

\newcommand{\bbb}[1]{\text{\boldmath$#1$}}
\newenvironment{pr}[1]{{\small\bf {#1}:}}{\hspace{\fill}$\square$}

\title{\Large\bf Moment-based inference for Pearson's quadratic \emph{q} subfamily of distributions}

\author{{\large G.\ AFENDRAS}\footnote{E-mail address:\ \tt g\_afendras@math.uoa.gr}}

\date{\small Department of Mathematics, Section of Statistics and O.R., University of Athens,\\
             Panepistemiopolis, 157 84 Athens, Greece.}

\begin{document}
 \thispagestyle{empty}
 \maketitle
 \vspace*{-2.5em}

 \begin{abstract}\vspace*{-.5ex}
 \noindent
 The author uses a Stein-type covariance identity to obtain moment estimators for the parameters of the quadratic polynomial subfamily of Pearson distributions. The asymptotic distribution of the estimators is obtained, and normality and symmetry tests based on it are provided. Simulation is used to compare the performance of the proposed tests with that of other existing tests for symmetry and normality.
 \end{abstract}
 \vspace*{-1ex}

 \noindent
 {\footnotesize {\sc Mathematics Subject Classification:} \ 62E01.
 }
 \vspace*{-1ex}

 \section{Introduction}\label{sec1}\vspace*{-.5em}
 Let $X$ be a continuous random variable (r.v.) with probability  density function (p.d.f.) $f$ and finite mean $\mu$. We say that $f$ has a Pearson quadratic form (see \cite{Kor}), $q(x)=\delta(x-\mu)^2+\beta(x-\mu)+\gamma$, if it satisfies the identity
 \begin{equation}\label{qua1}
 \int_{-\infty}^{x}(\mu-t)f(t)\d t=q(x)f(x) \ \ \ \text{for all} \ \ x\in\R.
 \end{equation}

 \begin{rem}\label{rem}
 Let $X$ be a continuous r.v. which satisfies {\rm(\ref{qua1})}. Then the support of $X$ is the interval
 \[
 J(X)=\big\{x: \ f(x)>0\big\}=\big(\essinf(X), \esssup(X)\big)=(\alpha,\omega),
 \]
 where $\essinf(X)=\inf_{x\in\R}\{F(x)>0\}$ and $\esssup(X)=\sup_{x\in\R}\{F(x)<1\}$. Also, it is obvious that $f\in{C}^\infty\big((\alpha,\omega)\big)$.
 \end{rem}

 \noindent
 The distributions satisfying (\ref{qua1}) belong to the Pearson family. Moreover, the quadratic $q$ generates the orthogonal polynomials through the Rodrigues-type formula, see \cite{Pap2},

 \[
 P_n(x)=\frac{(-1)^n}{f(x)}\frac{\d^n}{\d x^n}\big[q^n(x)f(x)\big], \ \ \ x\in J(X).
 \]

 Similarly, if $X$ is a discrete (integer-valued) r.v.\ with probability mass function (p.m.f.) $p$ and finite mean $\mu$, we say that $p$ has a Pearson quadratic form $q(j)=\delta(j-\mu)^2+\beta(j-\mu)+\gamma$ if it satisfies the identity
 \begin{equation}\label{qua2}
 \sum_{k\le{j}}(\mu-k)p(k)=q(j)p(j) \ \ \ \text{for all} \ \ j\in\Z.
 \end{equation}
 \begin{rem}\label{remm}
 Under {\rm(\ref{qua2})} it can be shown that the support \mbox{$J(X)=\big\{j\in\Z: p(j)>0\big\}$} is an interval of integers, i.e., if
 $j_1\in J(X)$ and $j_2\in J(X)$ then all integers $j$ between $j_1$ and $j_2$ belong to $J(X)$. The cases $J(X)=\{j_0\}$ or $J(X)=\{j_0,j_0+1\}$ are trivial, because identity {\rm(\ref{qua2})} is always satisfied and $q$ is not uniquely defined. We exclude these trivial cases from what follows. Thus, when we say that an integer-valued r.v.\ has a quadratic polynomial, we will assume that $\big|J(X)\big|\ge3$.
 \end{rem}

 For a suitable function, $g$, defined on $J(X)$, the following covariance identity holds (see \cite{APP2,CP1,Joh,Pap2}; cf.\ \cite{Ste1,Ste2}):
 \begin{equation}\label{c-iC}
 \E\big[(X-\mu)g(X)\big]=\E\big[q(X)g'(X)\big],\hspace{3ex}\text{(cont.\ case)}\vspace{-1ex}
 \end{equation}
 or\vspace{-1ex}
 \begin{equation}\label{c-iD}
 \E\big[(X-\mu)g(X)\big]=\E\big[q(X)\Delta{g(X)}\big],\hspace{3ex}\text{(discr.\ case)}
 \end{equation}
 where $\Delta$ denotes the forward difference operator, $\Delta g(x)=g(x+1)-g(x)$.

 Furthermore, using the Mohr and Noll inequality (or the discrete Mohr and Noll inequality), one obtains Poincar\'{e}-type lower/upper bounds for the variance of $g(X)$ (see \cite{APP1,Hou, HK,Joh,Pap1}) of the form
 \[
 (-1)^{n}\big[S_{n}-\Var g(X)\big]\geq0, \ \ \textrm{where}\vspace{-1ex}
 \]
 \[
 S_n=\sum_{k=0}^{n}\frac{(-1)^{k}}{(k+1)!\prod_{j=0}^{k}(1-j\delta)}\E\big[q^{k+1}(X)\left(g^{(k+1)}(X)\right)^2\big],
 \hspace{3ex}\text{(cont.\ case)} \vspace{-1ex}
 \]
 or\vspace{-1ex}
 \[
 S_n=\sum_{k=0}^{n}\frac{(-1)^{k}}{(k+1)!\prod_{j=0}^{k}(1-j\delta)}\E\big[q^{[k+1]}(X)\left(\Delta^{k+1}g(X)\right)^2\big],
 \hspace{3ex}\text{(discr.\ case)}
 \]
 where $\Delta^{k+1}g(x)=\Delta\left(\Delta^{k}g(x)\right)$, $k=1,2,\dots$, with $\Delta^0=I$.

 It should be noted that the quadratic $q$ also appears in variance bounds obtained using Bessel's inequality (see \cite{APP2}).

 Clearly, if we know the mean $\mu$ and the parameters $\delta,\beta,\gamma$ we can solve equation (\ref{qua1}) (or (\ref{qua2})) for $f$ (or $p$)
 (cf.\ \cite{CP1,Kor,Pap2,Pra1,Pra2}).

 The purpose of the present paper is to obtain an estimator for the parameters of the quadratic $q$, i.e., for the vector $(\mu,\delta,\beta,\gamma)^{\rm{t}}$. With the help of (\ref{c-iC}), (\ref{c-iD}), we generate a system of equations, from which we obtain the moment estimators for $\mu,\delta,\beta,\gamma$.

 Employing the $delta$-method, the asymptotic distribution of the estimators is derived. Some applications are also given. Similar
 work has been done by Pewsey, \cite{Pewsey}, who found the joint asymptotic distribution of the sample mean, variance, skewness and
 kurtosis. It is worth mentioning that Pewsey's results provide, for the first time, the joint asymptotic distribution for these fundamental statistics.

 \section{Moment Estimators}\vspace*{-.5em}
 \setcounter{equation}{0} \label{sec2}
 Here we deal with the estimation of the parameters $\delta,\beta,\gamma$ using the method of moments. ML estimation is possible but, as is generally true for all but the most simple of distributions, there are no closed-form expressions for the MLEs and ML estimation reduces to a numerical optimization problem. Instead, in what follows we consider estimators obtained using the method of moments.

 Let $X$  be a r.v.\ with $\E{X^4}<\infty$ and $\E{X}=\mu$. If $\mu_k=\E(X-\mu)^k$ is the $k$-th central moment, then
 \begin{equation}\label{l1}
 \mu_4\mu_2-\mu_3^2-\mu_2^3\ge0,
 \end{equation}
 and the equality holds only for the trivial case where $X$ takes, with probability $1$ (w.p.\ $1$), at most two values
 \big[the r.v.\ $Y=\mu_3(X-\mu)^2-(\mu_4-\mu_2^2)(X-\mu)$ has variance $\Var{Y}=(\mu_4-\mu_2^2)(\mu_4\mu_2-\mu_3^2-\mu_2^3)\ge0$, since $\Var(X-\mu)^2=\mu_4-\mu_2^2\ge0$\big]. Next, let $\bbb{X}_n=(X_1,\dots,X_n)$ be a random sample from any distribution. If  $m_{k;n}=\sum_{j=0}^{n}\left(X_j-\overline{X}_n\right)^k/n$ is the $k$-th sample central moment, then
 \begin{equation}\label{l2}
 m_{4;n}m_{2;n}-m_{3;n}^2-m_{2;n}^3\ge0,
 \end{equation}
 and the equality holds if and only if we have observed at most two values in the sample (this follows directly from (\ref{l1})).

\begin{theo}\label{theorem estimators}
 Let $X$ be an integer-valued r.v.\ (or a continuous r.v.) with mean $\mu$, finite fourth moment and p.m.f.\ $p$ satisfying  {\rm(\ref{qua2})} (or p.d.f.\ $f$ satisfying {\rm(\ref{qua1})}), with $q(j)=\delta(j-\mu)^2+\beta(j-\mu)+\gamma$. If $\bbb{X}_n=(X_1,\dots,X_n)$ is a random sample from $X$, with at least three different values (or at least three values), then
 \begin{itemize}
 \item[(a)] for the integer-valued case, the moment estimators for $\delta,\beta,\gamma$ are
 \[
 \begin{split}
 \widehat{\delta}_n&={\left(2m_{4;n}m_{2;n}-3m_{3;n}^2-6m_{2;n}^3+m_{2;n}^2\right)}\big/{6\widehat{\varTheta}_n},\\
  \widehat{\beta}_n&=\left[m_{3;n}\left(1-2\widehat{\delta}_n\right)-m_{2;n}\right]\big/2m_{2;n}\\
                   &={(m_{4;n}m_{3;n}-3m_{4;n}m_{2;n}+3m_{3;n}^2+3m_{3;n}m_{2;n}^2-m_{3;n}m_{2;n}+3m_{2;n}^3)}\big/{6\widehat{\varTheta}_n},\\
 \widehat{\gamma}_n&=m_{2;n}\left(1-\widehat{\delta}_n\right)={\left(4m_{4;n}m_{2;n}^2-3m_{3;n}^2m_{2;n}-m_{2;n}^3\right)}\big/{6\widehat{\varTheta}_n},
 \end{split}
 \]
 \item[(b)] for the continuous case, the moment estimators for $\delta,\beta,\gamma$ are
 \[
 \begin{split}
 \widehat{\delta}_n&={\left(2m_{4;n}m_{2;n}-3m_{3;n}^2-6m_{2;n}^3\right)}\big/{6\widehat{\varTheta}_n},\\
  \widehat{\beta}_n&=m_{3;n}\left(1-2\widehat{\delta}_n\right)\big/2m_{2;n}={\left(m_{4;n}m_{3;n}+3m_{3;n}m_{2;n}^2\right)}\big/{6\widehat{\varTheta}_n},\\
 \widehat{\gamma}_n&=m_{2;n}\left(1-\widehat{\delta}_n\right)={\left(4m_{4;n}m_{2;n}^2-3m_{3;n}^2m_{2;n}\right)}\big/{6\widehat{\varTheta}_n},
 \end{split}
 \]
 \end{itemize}
 where
 $\widehat{\varTheta}_n=m_{4;n}m_{2;n}-m_{3;n}^2-m_{2;n}^3$ is a positive number, $\overline{X}_{n}$ is the sample mean and $m_{k;n}=\sum_{j=1}^{n}\left(X_j-\overline{X}_n\right)^k/n$ is the $k$-th sample central moment. Also, the estimators $\widehat{\delta}_n,\widehat{\beta}_n,\widehat{\gamma}_n$ converge strongly to $\delta,\beta,\gamma$, respectively.
 \end{theo}

 \begin{pr}{Proof}
 (a) Since $X$ has finite fourth moment, it follows that it has finite central moments up to the fourth order. Also, its p.m.f.\ is
 quadratic $q$ and so the covariance identity (\ref{c-iD}) applies to any suitable $g$. In particular, for $g(x)=(x-\mu)^{k+1}, \ k=0,1,2$, the covariance identity is satisfied. For these functions, and since $\Delta(x-\mu)=1$, $\Delta(x-\mu)^2=2(x-\mu)+1$, $\Delta(x-\mu)^3=3(x-\mu)^2+3(x-\mu)+1$, (\ref{c-iD}) yields the equalities,
 \[
 \mu_2\delta+\gamma=\mu_2,
 \]
 \[
 (2\mu_3+\mu_2)\delta+2\mu_2\beta+\gamma=\mu_3,
 \]
 \[
 (3\mu_4+3\mu_3+\mu_2)\delta+(3\mu_3+3\mu_2)\beta+(3\mu_2+1)\gamma=\mu_4.
 \]
 Solving this system of equations, we obtain
 \[
 \begin{split}
 \delta&={\left(2\mu_4\mu_2-3\mu_3^2-6\mu_2^3+\mu_2^2\right)}\big/{6\varTheta},\\
  \beta&=[\mu_3(1-2\delta)-\mu_2]\big/2\mu_2
        ={\left(\mu_4\mu_3-3\mu_4\mu_2+3\mu_3^2+3\mu_3\mu_2^2-\mu_3\mu_2+3\mu_2^3\right)}\big/{6\varTheta},\\
 \gamma&=\mu_2(1-\delta)={\left(\mu_4\mu_2^2-\mu_3^2\mu_2-\mu_2^3\right)}\big/{6\varTheta},
 \end{split}
 \]
 where $\varTheta=\mu_4\mu_2-\mu_3^2-\mu_2^3$. For the solution of this system it is necessary to have $\varTheta\ne0$. This follows directly from (\ref{l1}), because $\big|J(X)\big|\ge3$ (see Remark \ref{remm}).

 \noindent
 If in a random sample we have observed at least three different values, then $\widehat{\varTheta}_n>0$, by (\ref{l2}). Replacing $m_{k;n}$ by $\mu_k$ in the above expressions we obtain the moment estimators $\widehat{\delta}_n$, $\widehat{\beta}_n$ and $\widehat{\gamma}_n$.
 \smallskip

 \noindent
 (b) Using similar arguments we observe that for $g(x)=(x-\mu)^{k+1}$, $k=0,1,2$, the covariance identity (\ref{c-iC}) is satisfied, which is
 \begin{equation}
 \mu_{k+2}=(k+1)\E\big[q(X)(X-\mu)^k\big]=(k+1)[\mu_{k+2}\delta+\mu_{k+1}\beta+\mu_{k}\gamma].\label{mu_k}
 \end{equation}
 From (\ref{mu_k}), with $k=0,1,2$, we generate a system of equations. Solving this system we obtain
 \[
 \begin{split}
 \delta&={\left(2\mu_2\mu_4-3\mu_3^2-6\mu_2^3\right)}\big/{6\varTheta},\\
  \beta&=\mu_3(1-2\delta)\big/2\mu_2={\left(\mu_4\mu_3+3\mu_3\mu_2^2\right)}\big/{6\varTheta},\\
 \gamma&=\mu_2(1-\delta)={\left(4\mu_4\mu_2^2-3\mu_3^2\mu_2\right)}\big/{6\varTheta},
 \end{split}
 \]
 where $\varTheta=\mu_4\mu_2-\mu_3^2-\mu_2^3$. For the solution of this system we have to assure that $\varTheta\ne0$. This follows directly from (\ref{l1}), since the r.v.\ $X$ is continuous.

 \noindent
 Since a random sample of at least three values from a continuous distribution function consists of distinct values (w.p.\ $1$), we have $\widehat{\varTheta}_n>0$, a.s.

 \noindent
 For all $k=1,2,3,4$ the r.v.\ $X^k$ has finite mean $\mu'_k=\E{X^k}$, and it is well known that  $m'_{k;n}=\sum_{j=1}^{n}{X_j^k}/n\longrightarrow\mu'_k$, a.s. Finally, the estimators $\widehat{\delta}_n,\widehat{\beta}_n,\widehat{\gamma}_n$ can  be written as rational functions of $m'_{k;n}, \ k=1,2,3,4$, and we conclude, using Slutsky's Theorem, that these functions converge strongly to $\delta,\beta,\gamma$ respectively.
 \end{pr}

 \begin{rem}\label{remak}
 If we carefully examine the expressions for $\delta,\beta,\gamma$ in the continuous case, we will see that $\delta$ is a number that does not have ``measurement units" (m.u.'s), $\beta$ is measured using the m.u.'s of $X$ and $\gamma$ is measured using the m.u.'s of the square of $X$. Bearing in mind that $\delta,\beta,\gamma$ are multiplied by $(x-\mu)^2,(x-\mu),1$ respectively, (i.e.\ $q(x)=\delta(x-\mu)^2+\beta(x-\mu)+\gamma$) we expect the final result to be measured in m.u.'s of the square of $X$. This indicates that the above choice of estimators is natural (see  {\rm(\ref{qua1})}).
 \end{rem}

 \section{Asymptotic Distribution}
 \setcounter{equation}{0}\label{sec3}
 Next, we study the asymptotic distribution of the estimators using the $delta$-method:
 \bigskip

 \noindent
 {\it Let $\bbb{T}_n$ be a sequence of r.v.'s in $\R^k$, $\bbb\vartheta\in\bbb\Theta\subseteq\R^k$, and assume that
 \[
 \sqrt{n}(\bbb{T}_n-\bbb{\vartheta})\lawv{N_k\left(\bbb0,\MATRIXS\right)}, \hspace{2ex} \text{as} \ \ n\to\infty,
 \]
 where $N_k(\bbb0,\MATRIXS)$ is a $k$-dimensional normal distribution with mean vector $\bbb0$ and covariance matrix $\MATRIXS$.

 \noindent
 If $\phi:\R^k\longrightarrow\R^m$ is (totally) differentiable in $\bbb\vartheta\in\bbb\Theta$, with total differential $J_\phi(\bbb{\vartheta})=\left.\left({\partial\phi_i(\bbb{x})}\big/{\partial{x}_j}\right)\right|_{\bbb{x}=\bbb{\vartheta}}\in\R^{m\times{k}}$, then (see {\rm\cite{Va}}, p.\ {\rm25})
 \[
 \sqrt{n}\big(\phi(\bbb{T}_n)-\phi(\bbb{\vartheta})\big)\lawv{N_m\left(\bbb0,J_\phi(\bbb{\vartheta})\MATRIXS J_\phi^{\rm{t}}(\bbb{\vartheta})\right)},
 \hspace{2ex} \text{as} \ \ n\to\infty.
 \]}
 So, we can easily deduce the following result.
 \begin{theo}\label{asymptotic}
 Let $\bbb{X}_n=(X_1,\dots,X_n)$ be a random sample, with at least three distinct values (or at least three values), from an integer-valued (or continuous) r.v.\ $X$, with p.m.f.\  satisfying {\rm(\ref{qua2})} (or p.d.f.\ satisfying {\rm(\ref{qua1})}) and finite $\E{X^8}$. If $\bbb{q}_n=\big(\overline{X}_n,\widehat{\delta}_n,\widehat{\beta}_n,\widehat{\gamma}_n\big)^{\rm{t}}$, with $\widehat{\delta}_n,\widehat{\beta}_n,\widehat{\gamma}_n$ as in Theorem {\rm\ref{theorem estimators}(a)} (or {\rm(b)}), and $\bbb{\omega}=(\mu,\delta,\beta,\gamma)^{\rm{t}}$, then

 \[
 \sqrt{n}(\bbb{q}_n-\bbb{\omega})\lawv{N_4(\bbb0,\mathbb{D})}, \hspace{2ex} \text{as} \ \ n\to\infty,
 \]

 \noindent
 where
 {\small$\mathbb{D}=\big(J_\varphi(\bbb{\vartheta})J_\psi(\bbb{\mu})\big)\MATRIXS\big(J_\varphi(\bbb{\vartheta})J_\psi(\bbb{\mu})\big)^{\rm{t}}$}\!, {\small$\bbb{\mu}=(0,\mu_2,\mu_3,\mu_4)^{\rm{t}}$}\!,
 $\bbb{\vartheta}=(\mu,\mu_2,\mu_3,\mu_4)^{\rm{t}}$,
 $\MATRIXS=(\mu_{i+j}-\mu_i\mu_j)_{i,j}$ (with $\mu_1=0$),
 $J_\psi(\bbb{\mu})$ has elements
 $j^\psi_{11}=1$,
  $j^\psi_{12}=0$,
   $j^\psi_{13}=0$,
    $j^\psi_{14}=0$,
 $j^\psi_{21}=0$,
  $j^\psi_{22}=1$,
   $j^\psi_{23}=0$,
    $j^\psi_{24}=0$,
 $j^\psi_{31}=-3\mu_2$,
  $j^\psi_{32}=0$,
   $j^\psi_{33}=1$,
    $j^\psi_{34}=0$,
 $j^\psi_{41}=-4\mu_3$,
  $j^\psi_{42}=0$,
   $j^\psi_{43}=0$,
    $j^\psi_{44}=1$,
 and $6\varTheta^2J_\varphi(\bbb{\vartheta})$ (with $\varTheta=\mu_4\mu_2-\mu_3^2-\mu_2^3$) has elements:
 \smallskip

 \noindent
 (a) for the discrete case,
 $j^\varphi_{11}=6\varTheta^2$,\,
  $j^\varphi_{12}=0$,\,
   $j^\varphi_{13}=0$,\,
    $j^\varphi_{14}=0$,\,
 $j^\varphi_{21}=0$,\,
  $j^\varphi_{22}=\mu_4\mu_3^2-8\mu_4\mu_2^3+\mu_4\mu_2^2+9\mu_3^2\mu_2^2-2\mu_3^2\mu_2+\mu_2^4$,\,
    $j^\varphi_{23}=-2\mu_4\mu_3\mu_2-6\mu_3\mu_2^3+2\mu_3\mu_2^2$,\,
     $j^\varphi_{24}=\mu_3^2\mu_2+4\mu_2^4-\mu_2^3$,\,
 $j^\varphi_{31}=0$,\,
  $j^\varphi_{32}=-\mu_4^2\mu_3+6\mu_4\mu_3\mu_2^2-6\mu_3^3\mu_2+\mu_3^3+3\mu_3\mu_2^4-\mu_3\mu_2^3$,\,
   $j^\varphi_{33}=\mu_4^2\mu_2+\mu_4\mu_3^2+2\mu_4\mu_2^3-\mu_4\mu_2^2+3\mu_3^2\mu_2^2-\mu_3^2\mu_2-3\mu_2^5+\mu_2^4$,\,
    $j^\varphi_{34}=-\mu_3^3-4\mu_3\mu_2^3+\mu_3\mu_2^2$,\,
 $j^\varphi_{41}=0$,\,
  $j^\varphi_{42}=\mu_4^2\mu_2^2-2\mu_4\mu_3^2\mu_2+\mu_4\mu_2^4-2\mu_4\mu_2^3+\mu_3^4-2\mu_3^2\mu_2^3+3\mu_3^2\mu_2^2$,\,
   $j^\varphi_{43}=2\mu_3\mu_2^4-2\mu_3\mu_2^3$,\,
    $j^\varphi_{44}=-\mu_2^5+\mu_2^4$,
 \smallskip

 \noindent
 (b) for the continuous case,
 $j^\varphi_{11}=6\varTheta^2$,\,\,
  $j^\varphi_{12}=0$,\,
   $j^\varphi_{13}=0$,\,
    $j^\varphi_{14}=0$,\,
 $j^\varphi_{21}=0$,\,
  $j^\varphi_{22}=\mu_4\mu_3^2-8\mu_4\mu_2^3+9\mu_3^2\mu_2^2$,\,
    $j^\varphi_{23}=-2\mu_4\mu_3\mu_2-6\mu_3\mu_2^3$,\,
     $j^\varphi_{24}=\mu_3^2\mu_2+4\mu_2^4$,\,
 $j^\varphi_{31}=0$,\,
  $j^\varphi_{32}=-\mu_4^2\mu_3+6\mu_4\mu_3\mu_2^2-6\mu_3^3\mu_2+3\mu_3\mu_2^4$,\,
   $j^\varphi_{33}=\mu_4^2\mu_2+\mu_4\mu_3^2+2\mu_4\mu_2^3+3\mu_3^2\mu_2^2-3\mu_2^5$,\,
    $j^\varphi_{34}=-\mu_3^3-4\mu_3\mu_2^3$,\,
 $j^\varphi_{41}=0$,\,
  $j^\varphi_{42}=4\mu_4^2\mu_2^2-8\mu_4\mu_3^2\mu_2+4\mu_4\mu_2^4+3\mu_3^4-6\mu_3^2\mu_2^3$,\,
   $j^\varphi_{43}=2\mu_4\mu_3\mu_2^2+6\mu_3\mu_2^4$,\,
    $j^\varphi_{44}=-\mu_3^2\mu_2^2-4\mu_2^5$.
 \end{theo}
 \begin{pr}{Proof}
 We centralize the $X_j$-values as $Y_j=X_j-\mu$, for $j=1,\dots,n$. Then for the vector $\bbb{Y}_n=\big(\overline{Y}_n,\overline{Y^2_n},\overline{Y^3_n},\overline{Y^4_n}\big)^{\rm{t}}$, it is well known that
 \[
 \sqrt{n}(\bbb{Y}_n-\bbb{\mu})\lawv N_4(\bbb0,\MATRIXS), \hspace{2ex} \text{as} \ \ n\to\infty.
 \]
 We consider the sample central moments $m_{k;n}=\sum_{j=1}^n\left(X_j-\overline{X}_n\right)^k\big/n$, $k=2,3,4$ and we seek the asymptotic distribution of the vector $\bbb{m}_n=\left(\overline{Y}_n,m_{2;n},m_{3;n},m_{4;n}\right)^{\rm{t}}$. Observe that
 \[
 m_{k;n}=\sum_{i=0}^{k}(-1)^{i}\binom{k}{i}\overline{Y_n^{k-i}}\left(\overline{Y}_n\right)^{i}.
 \]
 Hence, for $k=2,3,4$, we get
 $m_{2;n}=\overline{Y^2_n}-\left(\overline{Y}_n\right)^2$,
 $m_{3;n}=\overline{Y^3_n}-3\overline{Y^2_n}\left(\overline{Y}_n\right)+2\left(\overline{Y}_n\right)^3$,
 $m_{4;n}=\overline{Y^4_n}-4\overline{Y^3_n}\left(\overline{Y}_n\right)+6\overline{Y^2_n}\left(\overline{Y}_n\right)^2-3\left(\overline{Y}_n\right)^4$.
 Thus, the vector $\bbb{m}_n$ can be written as
 $\bbb{m}_n=\psi\big(\overline{Y}_n,\overline{Y^2_n},\overline{Y^3_n},\overline{Y^4_n}\big)$,
 where $\psi{(\bbb{x})}=\big(\psi_1(\bbb{x}),\psi_2(\bbb{x}),\psi_3(\bbb{x}),\psi_4(\bbb{x})\big)^{\rm{t}}$ with
 $\psi_1(\bbb{x})=x_1$,
  $\psi_2(\bbb{x})=x_2-x_1^2$,
   $\psi_3(\bbb{x})=x_3-3x_2x_1+2x_1^3$,
    $\psi_4(\bbb{x})=x_4-4x_3x_1+6x_2x_1^2-3x_1^4$.
 Applying the $delta$-method, it follows that
 $\sqrt{n}(\bbb{m}_n-\bbb{\mu})\lawv N_4\left(\bbb0,J_\psi(\bbb{\mu})\MATRIXS J_\psi^{\rm{t}}(\bbb{\mu})\right)$, as $n\to\infty$,
 where
 $J_\psi(\bbb{\mu})=\left({\partial\psi_i}\big/{\partial{x_j}}\right)\big|_{\bbb{x}=\bbb{\mu}}$.
 \smallskip

 \noindent
 Since $\bbb{T}_n=\left(\overline{X}_n,m_{2;n},m_{3;n},m_{4;n}\right)^{\rm{t}}=\bbb{m}_n+(\mu,0,0,0)^{\rm{t}}$ (and $\bbb{\vartheta}=\bbb{\mu}+(\mu,0,0,0)^{\rm{t}}$),
 we obtain
 $\sqrt{n}(\bbb{T}_n-\bbb{\vartheta})\lawv N_4\left(\bbb0,J_\psi(\bbb{\mu})\MATRIXS J_\psi^{\rm{t}}(\bbb{\mu})\right)$, as $n\to\infty$.
 \smallskip

 \noindent
 Regarding the asymptotic distribution of the vector $\sqrt{n}(\bbb{q}_n-\bbb{\omega})$,
 we have $\bbb{q}_n=\varphi(\bbb{T}_n)$ and $\bbb{\omega}=\varphi(\bbb{\vartheta})$,
 where $\varphi{(\bbb{x})}=\big(\varphi_1(\bbb{x}),\varphi_2(\bbb{x}),\varphi_3(\bbb{x}),\varphi_4(\bbb{x})\big)^{\rm{t}}$
 and the coordinates $\varphi_i(x)$ are given by:
 \begin{itemize}
 \item[(a)] for the integer-valued case,
 $\varphi_1(\bbb{x})=x_1$,
  $\varphi_2(\bbb{x})=(2x_4x_2-3x_3^2-6x_2^3+x_2^2)/c(\bbb{x})$,
   $\varphi_3(\bbb{x})=(x_4x_3-3x_4x_2+3x_3^2+3x_3x_2^2-x_3x_2+3x_2^3)/c(\bbb{x})$,
    $\varphi_4(\bbb{x})=(x_4x_2^2-x_3^2x_2-x_2^3)/c(\bbb{x})$,
 \item[(b)] for the continuous case,
 $\varphi_1(\bbb{x})=x_1$,
  $\varphi_2(\bbb{x})=\ds{\big(2x_4x_2-3x_3^2-6x_2^3\big)}/{c(\bbb{x})}$,
   $\varphi_3(\bbb{x})=\ds\big(x_4x_3$ $+3x_3x_2^2\big)/{c(\bbb{x})}$,
    $\varphi_4(\bbb{x})=\ds{\big(4x_4x_2^2-3x_3^2x_2\big)}/{c(\bbb{x})}$,
 \end{itemize}
 with $c(\bbb{x})=6(x_4x_2-x_3^2-x_2^3)$. Thus,
 \[
 \sqrt{n}(\bbb{q}_n-\bbb{\omega})\lawv
 N_4\left(\bbb0,\big(J_\varphi(\bbb{\vartheta})J_\psi(\bbb{\mu})\big)\MATRIXS J_\varphi(\bbb{\vartheta})J_\psi^{\rm{t}}(\bbb{\mu})\right),
 \hspace{2ex} \text{as} \ \ n\to\infty,
 \]
 where $J_\varphi(\bbb{\vartheta})=\left({\partial\varphi_i}\big/{\partial{x_j}}\right)\big|_{\bbb{x}=\bbb{\vartheta}}$ for both cases.
 \end{pr}

 \section{Hypothesis Testing}
 \setcounter{equation}{0}\label{sec4}
 In the subsections that follow we present various hypothesis tests based on the asymptotic distribution of the parameter estimates.

 \subsection{Continuous Case}\label{subsec4.1}

 \subsubsection{Test for Normality}\label{subsubsec4.1.1}
 A test of normality is equivalent to testing
 \[
 H_0: \delta=\beta=0 \ \ \text{vs} \ \
 H_1: \text{At least one of $\delta$ or $\beta$ is non-zero}.
 \]
 Theorem \ref{asymptotic}(b) shows that $\sqrt{n}\big(\widehat{\delta}_n-\delta,\widehat{\beta}_n-\beta\big)^{\rm{t}}\lawv{N_2}(\bbb0,\mathbb{D}_{\delta,\beta})$, as $n\to\infty$, where $\mathbb{D}_{\delta,\beta}=(AJ_\varphi(\bbb{\vartheta})J_\psi(\bbb{\mu}))\MATRIXS(AJ_\varphi(\bbb{\vartheta})J_\psi(\bbb{\mu}))^{\rm{t}}$, with
 $A=$ {\scriptsize$\Big(
   \begin{array}{@{\hspace{0ex}}c@{\hspace{1.1ex}}c@{\hspace{1.1ex}}c@{\hspace{1.1ex}}c@{\hspace{0ex}}}
   0 & 1 & 0 & 0 \\
   [-.8ex]
   0 & 0 & 1 & 0 \\
   \end{array}
   \Big)$}.
 Thus, under null hypothesis, we have that

 \[
 \ds{Q_n}=n\big(\widehat{\delta}_n,\widehat{\beta}_n,\big)\mathbb{D}^{-1}_{\delta,\beta;0}\big(\widehat{\delta}_n,\widehat{\beta}_n\big)^{\rm{t}}\lawv\chi^2_2,
 \hspace{2ex} \text{as} \ \ n\to\infty,
 \]

 \noindent
 where
 $\mathbb{D}_{\delta,\beta;0}$
  {\scriptsize
  $=\Big(
   \begin{array}{@{\hspace{0ex}}c@{\hspace{1.1ex}}c@{\hspace{0ex}}}
   2/3 & 0 \\
   [-.4ex]
   0 & 3\sigma^2/2 \\
   \end{array}
   \Big)$}
 and $\chi^2_m$ is the chi-square distribution with $m$ degrees of freedom. Since $\mathbb{D}_{\delta,\beta;0}$ is unknown we estimate it by $\widehat{\mathbb{D}_{\delta,\beta;0}}=$
 {\scriptsize
  $\Big(
   \begin{array}{@{\hspace{0ex}}c@{\hspace{1.1ex}}c@{\hspace{0ex}}}
   2/3 & 0 \\
   [-.4ex]
   0 & 3s^2/2 \\
   \end{array}
   \Big)$},
 replacing $\sigma^2$ by $s^2=\frac{n}{n-1}m_{2;n}$. For testing the above hypothesis, we propose the statistic
 \[
 q_n=n\big(\widehat{\delta}_n,\widehat{\beta}_n,\big)\widehat{\mathbb{D}_{\delta,\beta;0}}^{-1}\big(\widehat{\delta}_n,\widehat{\beta}_n\big)^{\rm{t}},
 \]
 and, at significance level $\alpha$, the asymptotic rejection region is $R=\left\{ q_n>\chi^2_{2;\alpha}\right\}$, where $\chi^2_{m;\alpha}$ is the upper $100\alpha\%$ point of the $\chi_m^2$ distribution.

 The distribution of $q_n$ is asymptotically $\chi^2_2$. Table 1 contains the 90th, 95th, 97.5th and 99th  percentiles of the empirical distribution  of $q_n$ generated by simulation of $10^5$ samples of size $n=10,20,30,50,70,100,150,200,300,400,500,750$ and $1000$ from a normal distribution.
 \smallskip

 \noindent
 {\footnotesize {\bf Table 1.} Empirical percentiles of the distribution of $q_n$ for random samples of size $n$ drawn from a normal distribution.}

 \centerline{\footnotesize
 \begin{tabular}{@{\hspace{0ex}}l@{\hspace{1.2ex}}c@{\hspace{1.2ex}}c@{\hspace{1.2ex}}c@{\hspace{1.2ex}}c@{\hspace{1.2ex}}c@{\hspace{1.2ex}}c@{\hspace{1.2ex}}c@{\hspace{1.2ex}}c@{\hspace{1.2ex}}c@{\hspace{1.2ex}}c@{\hspace{1.2ex}}c@{\hspace{1.2ex}}c@{\hspace{1.2ex}}c@{\hspace{1.2ex}}|c@{\hspace{0ex}}}
    \hline
    $n$ & $10$ & $20$ & $30$ & $50$ & $70$ & $100$ & $150$ & $200$ & $300$ & $400$ & $500$ & $750$ & $1000$ & $\infty$\\
    [.5ex]
   $P_{0.90}$ & $17.26$&$\ \ 9.18$&$\ \ 7.31$&$\ \ 6.06$&$\ \ 5.53$&$\ \ 5.20$&$\ \ 4.91$&$\ \ 4.85$&$4.75$&$4.66$&$4.67$&$4.63$&$4.62$&$4.60$\\
    [.5ex]
   $P_{0.95}$ & $27.48$&   $13.35$&   $10.37$&$\ \ 8.36$&$\ \ 7.56$&$\ \ 7.03$&$\ \ 6.53$&$\ \ 6.40$&$6.23$&$6.14$&$6.11$&$6.06$&$6.03$&$5.99$\\
   [.5ex]
   $P_{0.975}$& $41.25$&   $18.24$&   $13.83$&   $10.94$&$\ \ 9.75$&$\ \ 8.99$&$\ \ 8.23$&$\ \ 8.03$&$7.83$&$7.61$&$7.58$&$7.54$&$7.47$&$7.38$\\
   [.5ex]
   $P_{0.99}$ & $67.70$&   $26.03$&   $19.04$&   $14.52$&   $12.77$&   $11.78$&   $10.69$&   $10.33$&$9.96$&$9.72$&$9.55$&$9.42$&$9.39$&$9.21$\\
   \hline
 \end{tabular}}
 \medskip

 \noindent
 For a sample of small size $n$, the proposed $\alpha$-level normality test is
 \[
 \textrm{reject normality if \ } q_n\ge P_{1-\alpha},
 \]
 where $P_{1-\alpha}$ is given in Table 1.

 \subsubsection{Test for \bbb{\delta=0}}\label{subsubsec4.1.2}
 It is of interest to know if $\delta=0$, because this simplifies the procedure of inverting the quadratic $q$ and arranging (categorizing) the distribution. Hence we consider testing the null hypothesis $\delta=0$. Theorem \ref{asymptotic}(b) shows that $\sqrt{n}\big(\widehat{\delta}_n-\delta)\lawv{N}\big(0,\Var\widehat{\delta}_n\big)$, as $n\to\infty$, where $\Var\widehat{\delta}_n$ is the $(2,2)$ element of matrix $\mathbb{D}$.

 \noindent
 Note that if $q$ in (\ref{qua1}) is linear (that is, $\delta=0$) then $X$ follows either a normal or a gamma-type distribution of the form $X=cY+d$, where $Y$ is gamma and $c\ne0,d$ are constants. In both cases,
 $\Var\big(\widehat{\delta}_n\big|\delta=0\big)=({2}/{3})\big[1+{13\mu_3^2}\big/{4\mu_2^3}+{7\mu_3^4}\big/{2\mu_2^3(4\mu_2^3+\mu_3^2)}\big] \equiv\sigma_0^2\big(\widehat{\delta}_n\big)$.
 Thus, under null hypothesis,
 \[
 Z_n=\sqrt{n}\widehat{\delta}_n\big/\sigma_0\big(\widehat{\delta}_n\big)\lawv N(0,1), \hspace{2ex} \text{as} \ \ n\to\infty.
 \]
 However, $\sigma^2_0\big(\widehat{\delta}_n\big)$ is unknown, we have to estimate it by $\widehat{\sigma^2_0\big(\widehat{\delta}_n\big)}$, replacing $m_{k;n}$ by $\mu_k$. Thus, we proposed the statistic
 \[
 z_n=\sqrt{n}\frac{\widehat{\delta}_n}{\widehat{\sigma_0\big(\widehat{\delta}_n\big)}},
 \]
 and, at significance level $\alpha$, the (asymptotic) rejection region is $R=\left\{|z_n|>z_{\alpha/2}\right\}$, where $z_{\alpha}$ is the upper $100\alpha\%$ point of the standard normal distribution.

 \subsubsection{Test for Symmetry}\label{subsubsec4.1.3}
 First we prove the following lemma.
 \begin{lem}\label{lem symmetry c}
 Let $X$  be a continuous r.v.\ with mean $\mu$ and p.d.f.\ $f$ satisfying {\rm(\ref{qua1})}, with $q(x)=\delta(x-\mu)^2+\beta(x-\mu)+\gamma$. Then, $X$ is symmetric if and only if $\beta=0$.
 \end{lem}
 \begin{pr}{Proof}
 Let $X$ be symmetric. In this case $\mu=\E{X}$ is the point of symmetry. Taking derivatives in (\ref{qua1}) we get $(x-\mu)f(x)=q'(x)f(x)+q(x)f'(x)$ for all $x\in{J(X)}$. Since $\mu\in{J(X)}$, $(\mu-\mu)f(\mu)=q'(\mu)f(\mu)+q(\mu)f'(\mu)$. By the symmetry of $X$, $f'(\mu)=0$. So $q'(\mu)f(\mu)=0$ and since $f(x)>0$ for all $x\in{J(X)}$ we get $f(\mu)>0$ and $q'(\mu)=2\delta(\mu-\mu)+\beta=\beta=0$.

 \noindent
 Conversely, let $\beta=0$. Then $Y=X-\mu$ has $\mu_Y=\E{Y}=0$ and its density satisfies (\ref{qua1}), with $q_Y(y)=\delta{y^2}+\gamma$. We only have to show that $Y$ is symmetric at 0. The r.v.\ $-Y$ has $\mu_{-Y}=0$ and its density satisfies (\ref{qua1}), with $q_{-Y}(y)=q_Y(-y)=q_Y(y)$. \mbox{That is,} they have the same mean and the same Pearson quadratic form. Hence, $-Y\law{Y}$.
 \end{pr}
 \medskip

 We note that if $f$ satisfying (\ref{qua1}) where $\delta>0$ then $\E|X|^a<\infty$ for some $a>1$ if and only if $a<1+1/\delta$, see \cite{APP2}. Therefore, $\beta=0\Leftrightarrow\mu_3=0$ and a test of symmetry is equivalent to testing $H_0: \mu_3=0$ vs $H_1: \mu_3\ne0$. From the proof of Theorem \ref{asymptotic}(b) we have that $\sqrt{n}(m_{3;n}-\mu_3)^{\rm{t}}\lawv N(0,\sigma^2_{m_{3;n}})$, as $n\to\infty$, where $\sigma^2_{m_{3;n}}=\mu_6-6\mu_4\mu_2-\mu_3^2+9\mu_2^3$. Thus, under null hypothesis, $Z_n=\sqrt{n}m_{3;n}\big/\sigma_{m_{3;n};0}\lawv N(0,1)$, as $n\to\infty$, where (since $\mu_{2r+1}=0$) $\sigma^2_{m_{3;n};0}=\mu_6-6\mu_4\mu_2+9\mu_2^3$. Since $\sigma^2_{m_{3;n};0}$ is unknown, we estimate it by $\widehat{\sigma^2_{m_{3;n};0}}$, replacing $m_{2r;n}$ by $\mu_{2r}$, and the proposed statistic is
 \[
 z_n=\sqrt{n}\frac{m_{3;n}}{\widehat{\sigma_{m_{3;n};0}}},
 \]
 with asymptotic rejection region $R=\left\{|z_n|>z_{\alpha/2}\right\}$, at significance level $\alpha$.

 \subsection{Discrete Case}\label{subsec4.2}
 \subsubsection{Test for Poisson Distribution}\label{subsubsec4.2.1}
 An integer-valued r.v.\ $X$ with p.m.f.\ satisfying (\ref{qua2}), follows a Poisson distribution if and only if $\delta$, $\beta$ and $\sigma^2-\mu$ equal to zero.  Consider the test
 \[
 H_0: \delta=\beta=\sigma^2-\mu=0 \ \ \text{vs} \ \ H_1: \text{At least one of $\delta$ or $\beta$ or $\sigma^2-\mu$ is non-zero}.
 \]
 Let a function $\tau=(\tau_1,\tau_2,\tau_3)^{\rm{t}}:\R^4\longrightarrow\R^3$ with $\tau_i(\bbb{x})\equiv\varphi_{i+1}(\bbb{x})$, $i=1,2$ (where  $\varphi_{i+1}$ as in proof of Theorem \ref{asymptotic}(a)) and $\tau_3(\bbb{x})=x_2-x_1$.
 Note that $\tau(\bbb{\vartheta})=(\delta,\beta,\sigma^2-\mu)^{\rm{t}}$ and $\tau(\bbb{T}_n)=\big(\widehat{\delta}_n,\widehat{\beta}_n,m_{2;n}-\overline{X}_n\big)^{\rm{t}}$,
 where $\widehat{\delta}_n$, $\widehat{\beta}_n$ the estimators in Theorem \ref{theorem estimators}(a).
 Therefore, a test for Poisson is reduced to the null hypothesis
 $\tau(\bbb{\vartheta})=\bbb{0}$. Using the $delta$-method, $\sqrt{n}\big(\tau\big(\widehat{\bbb{T}}_n\big)-\tau(\bbb{\vartheta})\big)\lawv N_3\big(0,\mathbb{D}_\tau\big)$, as $n\to\infty$,
 where $\mathbb{D}_\tau=C(\bbb{\vartheta})\MATRIXS{C^{\rm{t}}(\bbb{\vartheta})}$, with $C(\bbb{\vartheta})=J_\tau(\bbb{\vartheta})J_\psi(\bbb{\mu})$
 and $J_\tau(\bbb{\vartheta})=\left({\partial\tau_i}\big/{\partial{x_j}}\right)\big|_{\bbb{x}=\bbb{\vartheta}}$.
 The first two rows of matrix $J_\tau(\bbb{\vartheta})$ are the rows two and three of matrix $J_\varphi(\bbb{\vartheta})$ in Theorem \ref{asymptotic}(a), the third row of this matrix is $(-1,1,0,0)$. Thus, under null hypothesis,
 \[
 Q_n=n\big(\widehat{\delta}_n,\widehat{\beta}_n,m_{2;n}-\overline{X}_n\big)\mathbb{D}^{-1}_{\tau;0} \big(\widehat{\delta}_n,\widehat{\beta}_n,m_{2;n}-\overline{X}_n\big)^{\rm{t}}
 \lawv\chi^2_3,\hspace{2ex} \text{as} \ \ n\to\infty.
 \]
 where
 $\mathbb{D}_{\tau;0}=\frac{1}{6\lambda}${\scriptsize$\bigg(
   \begin{array}{@{\hspace{0ex}}c@{\hspace{2ex}}c@{\hspace{2ex}}c@{\hspace{0ex}}}
   4\lambda+9 & 5\lambda-9 & 0 \\
   [-.6ex]
   5\lambda-9 & 9\lambda^2-2\lambda & 12\lambda^2 \\
   [-.6ex]
   0 & 12\lambda^2 & 12\lambda^3 \\
   \end{array}
   \bigg)$}
 and $\lambda$ is the parameter of the Poisson distribution. Since $\mathbb{D}_{\tau;0}$ is unknown, we estimate it by $\widehat{\mathbb{D}_{\tau;0}}$, replacing $\lambda$ by $\overline{X}_n$, and the proposed statistic is
 \[
 q_n=n\big(\widehat{\delta}_n,\widehat{\beta}_n,m_{2;n}-\overline{X}_n\big)\widehat{\mathbb{D}_{\tau;0}}^{-1}
 \big(\widehat{\delta}_n,\widehat{\beta}_n,m_{2;n}-\overline{X}_n\big)^{\rm{t}},
 \]
 with asymptotic rejection region $R=\left\{q_n>\chi^2_{3;\alpha}\right\}$, at significance level $\alpha$.
 \bigskip

 \section{Simulation Results}\label{sec5}
 {\small\bf Continuous case} \
 Using Matlab, we simulated $10^4$ samples of size $n$ from various continuous distributions admitting a Pearson quadratic form (see Figure 1).

 In Tables 2--10 we present the observed averages of the estimators $\widehat{\delta}_n,\widehat{\beta}_n$ and $\widehat{\gamma}_n$ (in Theorem \ref{theorem estimators}(b)), their observed mean square errors ($MSEs$) and the empirical level of rejection for the test of $\delta=0$ for a nominal significance level of $\alpha=0.05$. Considering the content of these Tables, we see that, as $n$ becomes large, the means of the estimators tend to their true values. It would appear that the estimators of $\delta$ and $\gamma$ have negative and positive  bias, respectively. Also, as expected, the mean square error of these estimators tends to zero, as $n$ tends to infinity.

 There exists a well-established literature addressing the problem of testing univariate data for normality. Renewed recent interest in this inferential problem can be found in \cite{TB,ZW}, see also \cite{JP}. In the light of the findings presented in those papers, we conducted a simulation study designed to compare the performance of the new proposed test of normality with those of the following competitive tests for a nominal significance level of 5\%, see Figures 2(a)--10(a). The description of each test is preceded by the abbreviation we will use when referring to it:
 \begin{itemize}
 \item[]
 \begin{itemize}
 \item[KS.] The (one-sided) empirical distribution function (e.d.f.)-based test of  Kolmo\-gorov--Smirnov.
 \item[BS.] The test of Bowman and Shenton (see \cite{BS}), or the (one-sided) Jarque and Bera (see \cite{JB}), the test statistic of which is a function of the coefficients of skewness and kurtosis.
 \item[D.]  The (two-sided) test of D' Agostino (see \cite{DAg1,DAg2}). Up to a constant, the test statistic is the ratio of Downton's (see \cite{Dow}) linear estimator of the standard deviation to the sample standard deviation. The critical values for this test are given in D' Agostino's papers.
 \item[AD.] The (one-sided) empirical distribution function (e.d.f.)-based test of Anderson and Darling (see \cite{And-Darl}). We used the corrected critical values for this test presented under the name CMWS in Table 2 of see \cite{TB}.
 \item[CvM.]The (one-sided) Cram\'{e}r-von Mises e.d.f.-based test with statistic identified as CMS in \cite{TB}. We used the corrected  critical values given in their Table 2.
 \item[$Z_A$, $Z_C$.] The (one-sided) nonparametric likelihood-ratio-based tests with test statistics $Z_A$ and $Z_C$ of Zhang and Wu (see \cite{ZW}). We used the corrected critical values for these tests given in Tables 1 and 2, respectively, of that paper.
 \end{itemize}
 \end{itemize}

 We will compare the size and power performance of the proposed symmetry test (again, for a nominal significance level of 0.05) with that of one other general test of symmetry, see Figures 2(b)--10(b). This particular test were chosen because they were found to perform well in extensive simulation comparisons reported in Cabilio and Masaro (1996). That is:
 \begin{itemize}
 \item[CM.] The test statistic of Cabilio and Masaro (see \cite{Cab-Mas}), is the simple function  $S_K=\sqrt{n}(\overline{X}-m)\big/s$, where $\overline{X}$, $m$ and $s$ denote the sample mean, median and standard deviation (with divisor $n$),  respectively. Under symmetry, $S_K\lawv N\big(0,\sigma_0^2(F)\big)$, as $n\to\infty$, where $\sigma_0^2(F)$ is a constant (depending on the distribution $F$). The  critical values for this test are given in \cite{Cab-Mas}. Notice that this dependence is a weak point of the test statistic.
 \end{itemize}
 \vspace{2cm}

 Figure 1 portrays the densities of the variates used to produce the numerical results presented in Tables 2--10 and Figures 2--10.
 \smallskip

 \centerline{\includegraphics[scale=1]{Figure.1}}

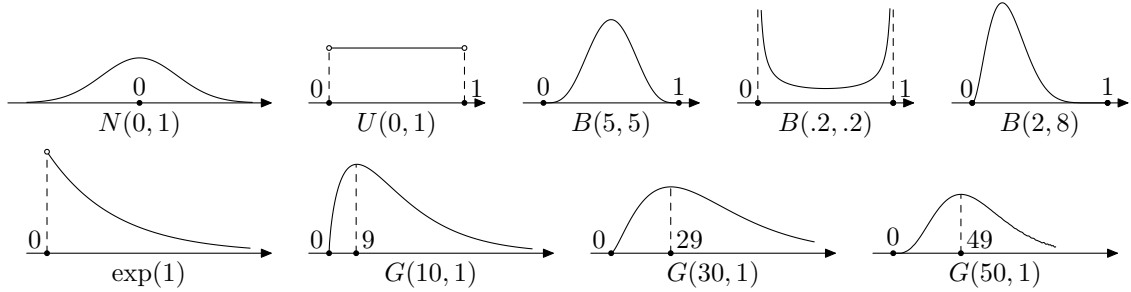
 \centerline{\footnotesize {\bf Figure 1.} Densities of standard normal and various beta and gamma distributions.}
 \bigskip

 For the normal distribution the parameters $\mu$ and $\sigma^2$ do not affect the critical region of the tests, so we fix $\mu=0$ and $\sigma^2=1$.
 \medskip

 \noindent
 {\footnotesize {\bf Table 2.} Mean and mean square error (MSE) for the three estimators calculated from $10^4$ random samples of size $n$ simulated from the standard normal distribution (i.e.\ $\delta=0$, $\beta=0$, $\gamma=1$). The last row gives the empirical size, $p_\delta$, of the test for $\delta=0$ with a nominal significance level of $\alpha=0.05$.}

 \centerline{\footnotesize
 \begin{tabular}{@{\hspace{.5ex}}l
                 @{\hspace{1.3ex}}r
                 @{\hspace{1.3ex}}r
                 @{\hspace{1.3ex}}r
                 @{\hspace{1.3ex}}r
                 @{\hspace{1.3ex}}r
                 @{\hspace{1.3ex}}r
                 @{\hspace{1.3ex}}r
                 @{\hspace{.5ex}}}
   \hline
          $n$                       &
          \multicolumn{1}{c}{$50$}  &
          \multicolumn{1}{c}{$100$} &
          \multicolumn{1}{c}{$150$} &
          \multicolumn{1}{c}{$200$} &
          \multicolumn{1}{c}{$300$} &
          \multicolumn{1}{c}{$400$} &
          \multicolumn{1}{c}{$500$} \\ 
          $\widehat{\delta}_n$      & $-0.0818$ & $-0.0417$ & $-0.0283$ & $-0.0208$ & $-0.0144$ & $-0.0114$ & $-0.0089$  \\
   [-.2ex]$MSE(\widehat{\delta}_n)$ & $ 0.0188$ & $ 0.0077$ & $ 0.0049$ & $ 0.0035$ & $ 0.0023$ & $ 0.0018$ & $ 0.0014$  \\
          $\widehat{\beta}_n$       & $-0.0018$ & $ 0.0015$ & $-0.0005$ & $-0.0002$ & $ 0.0013$ & $-0.0005$ & $-0.0000$  \\
   [-.2ex]$MSE(\widehat{\beta}_n)$  & $ 0.0336$ & $ 0.0160$ & $ 0.0103$ & $ 0.0077$ & $ 0.0051$ & $ 0.0038$ & $ 0.0030$  \\
          $\widehat{\gamma}_n$      & $ 1.0611$ & $ 1.0290$ & $ 1.0199$ & $ 1.0169$ & $ 1.0106$ & $ 1.0099$ & $ 1.0070$  \\
   [-.2ex]$MSE(\widehat{\gamma}_n)$ & $ 0.0614$ & $ 0.0281$ & $ 0.0181$ & $ 0.0139$ & $ 0.0091$ & $ 0.0069$ & $ 0.0054$  \\
          $p_{\delta}$              & $ 0.0778$ & $ 0.0616$ & $ 0.0577$ & $ 0.0516$ & $ 0.0516$ & $ 0.0476$ & $ 0.0504$  \\
   \hline
 \end{tabular}}
 \bigskip

 \centerline{\includegraphics[scale=1]{Figure.2}}

 \noindent
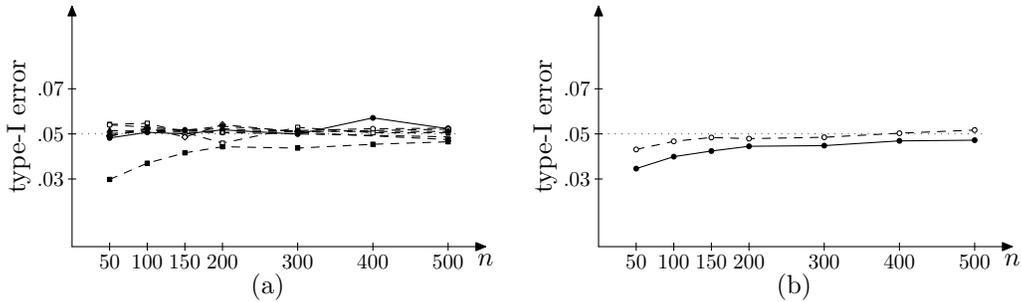
 {\footnotesize {\bf Figure 2.} Empirical type I error calculated from $10^4$ random samples of size $n$ simulated from the standard normal distribution, for (a) tests for normality and (b) tests for symmetry with a nominal size of 0.05.
   The results in (a) correspond to the proposed new test (\hspace{-.2ex}\CircleB) and the seven existing tests
   KS $\langle$\hspace{-.2ex}\Square$\rangle$,
   BS $\langle$\hspace{-.2ex}\SquareB$\rangle$,
   D  $\langle$\hspace{-.2ex}\Circle$\rangle$,
   AD $\langle$\Triangle$\rangle$,
   CvM $\langle$\TriangleB$\rangle$,
   $Z_A$ $\langle$\Diamond$\rangle$,
   $Z_C$ $\langle$\DiamondB$\rangle$.
   The results in (b) correspond to the proposed new test $\langle$\hspace{-.2ex}\CircleB$\rangle$, and the CM test $\langle$\hspace{-.2ex}\Circle$\rangle$.}
 \bigskip

 \noindent
 From Table 2, the test for $\delta=0$ maintains the nominal level well for a sample size of $n\ge200$ and is liberal for $n\le150$.
 \smallskip

 \noindent
 From Figure 2, the tests for normality hold the nominal level well, apart from the BS test which is conservative. As the sample size $n$ increases the empirical sizes of the tests for symmetry tend to the nominal level of $\alpha=0.05$, as expected. The CM test for symmetry maintains the nominal level better than the new proposed test for symmetry which is conservative.
 \smallskip

 Next we simulate data from beta and gamma distributions; clearly, they satisfy {\rm(\ref{qua1})}.
 \smallskip

 Let $X$ be beta distribution with parameters $a$ and $b$. The r.v.\ $X$ is symmetric if and only if $a=b$. Also, for some values of $a$ and $b$ the density function of $X$ is close to some normal distribution, for example see Figure 4.1. For various values of $a$, $b$ we have:
 \smallskip

 \noindent
 {\footnotesize {\bf Table 3.} Mean and mean square error (MSE) for the three estimators calculated from $10^4$ random samples of size $n$ simulated from  the standard uniform distribution (i.e.\ $\delta=-0.5$, $\beta=0$, $\gamma=0.125$). The last row gives the empirical power, $p^*_\delta$, of the test for $\delta=0$ with a nominal significance level of $\alpha=0.05$.}

 \centerline{\footnotesize
 \begin{tabular}{@{\hspace{.5ex}}l
                 @{\hspace{1.3ex}}r
                 @{\hspace{1.3ex}}r
                 @{\hspace{1.3ex}}r
                 @{\hspace{1.3ex}}r
                 @{\hspace{1.3ex}}r
                 @{\hspace{1.3ex}}r
                 @{\hspace{1.3ex}}r
                 @{\hspace{.5ex}}}
   \hline
          $n$                       &
          \multicolumn{1}{c}{$50$}  &
          \multicolumn{1}{c}{$100$} &
          \multicolumn{1}{c}{$150$} &
          \multicolumn{1}{c}{$200$} &
          \multicolumn{1}{c}{$300$} &
          \multicolumn{1}{c}{$400$} &
          \multicolumn{1}{c}{$500$} \\ 
          $\widehat{\delta}_n$      & $-0.5365$ & $-0.5170$ & $-0.5114$ & $-0.5093$ & $-0.5054$ & $-0.5039$ & $-0.5037$  \\
   [-.2ex]$MSE(\widehat{\delta}_n)$ & $ 0.0349$ & $ 0.0153$ & $ 0.0101$ & $ 0.0078$ & $ 0.0049$ & $ 0.0036$ & $ 0.0029$  \\
          $\widehat{\beta}_n$       & $ 0.0001$ & $-0.0004$ & $-0.0002$ & $ 0.0006$ & $-0.0001$ & $-0.0003$ & $-0.0001$  \\
   [-.2ex]$MSE(\widehat{\beta}_n)$  & $ 0.0038$ & $ 0.0018$ & $ 0.0012$ & $ 0.0009$ & $ 0.0006$ & $ 0.0004$ & $ 0.0003$  \\
          $\widehat{\gamma}_n$      & $ 0.1271$ & $ 0.1259$ & $ 0.1256$ & $ 0.1256$ & $ 0.1253$ & $ 0.1252$ & $ 0.1252$  \\
   [-.2ex]$MSE(\widehat{\gamma}_n)$ & $ 0.0009$ & $ 0.0004$ & $ 0.0003$ & $ 0.0002$ & $ 0.0001$ & $ 0.0001$ & $ 0.0001$  \\
          $p_{\delta}^*$            & $ 0.9712$ & $ 0.9997$ & $ 1.0000$ & $ 1.0000$ & $ 1.0000$ & $ 1.0000$ & $ 1.0000$  \\
   \hline
 \end{tabular}}
 \medskip

 \centerline{\includegraphics[scale=1]{Figure.3}}\vspace{-1.5ex}

 \noindent
 {\footnotesize {\bf Figure 3.} Empirical power of the tests for normality (a) and the size of tests for symmetry (b), for a nominal size of $0.05$, calculated from $10^4$ random samples of size $n$ simulated from the standard uniform distribution. The tests are the same as those described in the caption to Figure 2.}
 \smallskip

 \noindent
 From Table 3, as the sample size $n$ increases, the power of the test for $\delta=0$ increases to one rapidly.
 \smallskip

 \noindent
 From Figure 3, as the sample size  $n$ increases, the power of the tests for normality increases to one and the empirical size of the tests for symmetry tend to the nominal level of $\alpha=0.05$. For $n=50$ the proposed test for normality is more powerful than the other tests for normality. The tests BS (for a sample size of  $n=50$) and KS (for a sample size of $n\le300$) have a poor performance. The new proposed test for symmetry maintains the nominal level better than the CM test for symmetry which is conservative.
 \vspace{2cm}

 \noindent
 {\footnotesize {\bf Table 4.} Mean and mean square error (MSE) for the three estimators calculated from $10^4$ random samples of size $n$ simulated from the beta distribution with $a=b=5$ (i.e.\ $\delta=-0.1$, $\beta=0$, $\gamma=0.025$). The last row gives the empirical power, $p^*_\delta$, of the test for $\delta=0$ with a nominal significance level of $\alpha=0.05$.}

 \centerline{\footnotesize
 \begin{tabular}{@{\hspace{.5ex}}l
                 @{\hspace{1.3ex}}r
                 @{\hspace{1.3ex}}r
                 @{\hspace{1.3ex}}r
                 @{\hspace{1.3ex}}r
                 @{\hspace{1.3ex}}r
                 @{\hspace{1.3ex}}r
                 @{\hspace{1.3ex}}r
                 @{\hspace{.5ex}}}
   \hline
          $n$                       &
          \multicolumn{1}{c}{$50$}  &
          \multicolumn{1}{c}{$100$} &
          \multicolumn{1}{c}{$150$} &
          \multicolumn{1}{c}{$200$} &
          \multicolumn{1}{c}{$300$} &
          \multicolumn{1}{c}{$400$} &
          \multicolumn{1}{c}{$500$} \\ 
          $\widehat{\delta}_n$      & $-0.1599$ & $-0.1274$ & $-0.1190$ & $-0.1141$ & $-0.1086$ & $-0.1060$ & $-0.1052$  \\
   [-.2ex]$MSE(\widehat{\delta}_n)$ & $ 0.0169$ & $ 0.0066$ & $ 0.0041$ & $ 0.0030$ & $ 0.0019$ & $ 0.0014$ & $ 0.0011$  \\
          $\widehat{\beta}_n$       & $-0.0001$ & $-0.0000$ & $ 0.0001$ & $ 0.0001$ & $-0.0001$ & $ 0.0000$ & $-0.0000$  \\
   [-.2ex]$MSE(\widehat{\beta}_n)$  & $ 0.0006$ & $ 0.0003$ & $ 0.0002$ & $ 0.0001$ & $ 0.0001$ & $ 0.0001$ & $ 0.0001$  \\
          $\widehat{\gamma}_n$      & $ 0.0259$ & $ 0.0254$ & $ 0.0253$ & $ 0.0253$ & $ 0.0251$ & $ 0.0251$ & $ 0.0251$  \\
   [-.2ex]$MSE(\widehat{\gamma}_n)$ & $ 0.0000$ & $ 0.0000$ & $ 0.0000$ & $ 0.0000$ & $ 0.0000$ & $ 0.0000$ & $ 0.0000$  \\
          $p_{\delta}^*$            & $ 0.2073$ & $ 0.2812$ & $ 0.3762$ & $ 0.4574$ & $ 0.6193$ & $ 0.7560$ & $ 0.8368$  \\
   \hline
 \end{tabular}}
 \medskip

 \centerline{\includegraphics[scale=1]{Figure.4}}\vspace{-1.5ex}

 \noindent
 {\footnotesize {\bf Figure 4.} Empirical power of the tests for normality (a) and the size of tests for symmetry (b), for a nominal size of 0.05, calculated from $10^4$ random samples of size $n$ simulated from the beta distribution with $a=b=5$. The tests are the same as those described in the caption to Figure 2.}
 \medskip

 \noindent
 From Table 4, as the sample size  $n$ increases, the power of the test for $\delta=0$ increases. That test for a sample size of $n\le100$ has a poor performance.
 \smallskip

 \noindent
 From Figure 4, as the sample size  $n$ increases, the power of the tests for normality increases. The new proposed and D tests for normality are more powerful than the other tests for normality. All tests (new proposed for normality and D for a sample size of $n\le200$; $Z_A$ and $Z_C$ for a sample size of $n\le300$; BS, KS and AD for a sample size of $n\le500$) have a poor performance. This happens because the density of beta distribution $B(5,5)$ is close to the normal distribution $N(1/2,1/44)$ (see Figure 4.1). Notice that the power of the KS test is less than the nominal level for a sample size of $n\le500$. Both tests for symmetry maintain the nominal level equally well. Both are conservative.
 \medskip

 \centerline{\includegraphics[scale=1]{Figure.41}}
 \noindent
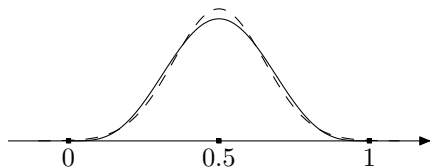
 {\footnotesize {\bf Figure 4.1.} Densities of the beta distribution with parameters $\alpha=5$, $\beta=5$ (solid line) and the normal distribution with parameters $\mu=1/2$, $\sigma^2=1/44$ (dashed line).}
 \vspace{2cm}

 \noindent
 {\footnotesize {\bf Table 5.} Mean and mean square error (MSE) for the three estimators calculated from $10^4$ random samples of size $n$ simulated from the beta distribution with $a=b=0.2$ (i.e.\ $\delta=-2.5$, $\beta=0$, $\gamma=0.625$). The last row gives the empirical power, $p^*_\delta$, of the test for $\delta=0$ with a nominal significance level of $\alpha=0.05$.}

 \centerline{\footnotesize
 \begin{tabular}{@{\hspace{.5ex}}l
                 @{\hspace{1.3ex}}r
                 @{\hspace{1.3ex}}r
                 @{\hspace{1.3ex}}r
                 @{\hspace{1.3ex}}r
                 @{\hspace{1.3ex}}r
                 @{\hspace{1.3ex}}r
                 @{\hspace{1.3ex}}r
                 @{\hspace{.5ex}}}
   \hline
          $n$                       &
          \multicolumn{1}{c}{$50$}  &
          \multicolumn{1}{c}{$100$} &
          \multicolumn{1}{c}{$150$} &
          \multicolumn{1}{c}{$200$} &
          \multicolumn{1}{c}{$300$} &
          \multicolumn{1}{c}{$400$} &
          \multicolumn{1}{c}{$500$} \\ 
          $\widehat{\delta}_n$      & $-2.6784$ & $-2.5798$ & $-2.5523$ & $-2.5310$ & $-2.5209$ & $-2.5220$ & $-2.5177$  \\
   [-.2ex]$MSE(\widehat{\delta}_n)$ & $ 0.9962$ & $ 0.3920$ & $ 0.2377$ & $ 0.1693$ & $ 0.1097$ & $ 0.0861$ & $ 0.0678$  \\
          $\widehat{\beta}_n$       & $ 0.0009$ & $-0.0010$ & $ 0.0005$ & $-0.0001$ & $ 0.0010$ & $-0.0015$ & $ 0.0003$  \\
   [-.2ex]$MSE(\widehat{\beta}_n)$  & $ 0.1266$ & $ 0.0568$ & $ 0.0358$ & $ 0.0261$ & $ 0.0172$ & $ 0.0129$ & $ 0.0104$  \\
          $\widehat{\gamma}_n$      & $ 0.6545$ & $ 0.6376$ & $ 0.6332$ & $ 0.6293$ & $ 0.6280$ & $ 0.6287$ & $ 0.6280$  \\
   [-.2ex]$MSE(\widehat{\gamma}_n)$ & $ 0.0505$ & $ 0.0203$ & $ 0.0124$ & $ 0.0088$ & $ 0.0058$ & $ 0.0045$ & $ 0.0035$  \\
          $p_{\delta}^*$            & $ 1.0000$ & $ 1.0000$ & $ 1.0000$ & $ 1.0000$ & $ 1.0000$ & $ 1.0000$ & $ 1.0000$  \\
   \hline
 \end{tabular}}
 \medskip

 \centerline{\includegraphics[scale=.99]{Figure.5}}\vspace{-.7ex}

 \noindent
 {\footnotesize {\bf Figure 5.} Empirical power of the tests for normality (a) and the size of tests for symmetry (b), for a nominal size of 0.05, calculated from $10^4$ random samples of size $n$ simulated from the beta distribution with $a=b=0.2$. The tests are the same as those described in the caption to Figure 2.}
 \medskip

 \noindent
 From Table 5, for each sample size $n$ the power of the test for $\delta=0$ is one.
 \smallskip

 \noindent
 From Figure 5, as the sample size  $n$ increases, the power of the tests for normality increases to one and the empirical size of the tests for symmetry tend to the nominal level of $\alpha=0.05$. The test D for a sample size of $n\le150$ has a poor performance. The new proposed test for symmetry maintains the nominal level better than the CM test for symmetry. However, the new test is liberal whereas the CM test is conservative.
 \medskip

 \noindent
 {\footnotesize {\bf Table 6.} Mean and mean square error (MSE) for the three estimators calculated from $10^4$ random samples of size $n$ simulated from the beta distribution with $a=2, b=8$ (i.e.\ $\delta=-0.1$, $\beta=0.06$, $\gamma=0.016$). The last row gives the empirical power, $p^*_\delta$, of the test for $\delta=0$ with a nominal significance level of $\alpha=0.05$.}

 \centerline{\footnotesize
 \begin{tabular}{@{\hspace{.5ex}}l
                 @{\hspace{1.3ex}}r
                 @{\hspace{1.3ex}}r
                 @{\hspace{1.3ex}}r
                 @{\hspace{1.3ex}}r
                 @{\hspace{1.3ex}}r
                 @{\hspace{1.3ex}}r
                 @{\hspace{1.3ex}}r
                 @{\hspace{.5ex}}}
   \hline
          $n$                       &
          \multicolumn{1}{c}{$50$}  &
          \multicolumn{1}{c}{$100$} &
          \multicolumn{1}{c}{$150$} &
          \multicolumn{1}{c}{$200$} &
          \multicolumn{1}{c}{$300$} &
          \multicolumn{1}{c}{$400$} &
          \multicolumn{1}{c}{$500$} \\ 
          $\widehat{\delta}_n$      & $-0.2052$ & $-0.1552$ & $-0.1385$ & $-0.1284$ & $-0.1204$ & $-0.1158$ & $-0.1116$  \\
   [-.2ex]$MSE(\widehat{\delta}_n)$ & $ 0.0298$ & $ 0.0128$ & $ 0.0082$ & $ 0.0059$ & $ 0.0040$ & $ 0.0031$ & $ 0.0024$  \\
          $\widehat{\beta}_n$       & $ 0.0622$ & $ 0.0617$ & $ 0.0610$ & $ 0.0609$ & $ 0.0606$ & $ 0.0604$ & $ 0.0603$  \\
   [-.2ex]$MSE(\widehat{\beta}_n)$  & $ 0.0009$ & $ 0.0004$ & $ 0.0002$ & $ 0.0002$ & $ 0.0001$ & $ 0.0001$ & $ 0.0001$  \\
          $\widehat{\gamma}_n$      & $ 0.0172$ & $ 0.0167$ & $ 0.0165$ & $ 0.0164$ & $ 0.0162$ & $ 0.0162$ & $ 0.0161$  \\
   [-.2ex]$MSE(\widehat{\gamma}_n)$ & $ 0.0000$ & $ 0.0000$ & $ 0.0000$ & $ 0.0000$ & $ 0.0000$ & $ 0.0000$ & $ 0.0000$  \\
          $p_{\delta}^*$            & $ 0.1511$ & $ 0.1606$ & $ 0.1678$ & $ 0.1909$ & $ 0.2403$ & $ 0.2943$ & $ 0.3520$  \\
   \hline
 \end{tabular}}
 \medskip

 \noindent
 From Table 6, as the sample size  $n$ increases, the power of the test for $\delta=0$ increases. That test for each sample size of $n$ has a poor performance.
 \medskip

 \centerline{\includegraphics[scale=1]{Figure.6}}

  \noindent
 {\footnotesize {\bf Figure 6.} Empirical power of the tests for normality (a) and symmetry (b), for a nominal size of 0.05, calculated from $10^4$ random samples of size $n$ simulated from the beta distribution with $a=2, b=8$. The tests are the same as those described in the caption to Figure 2.}
 \medskip

 \noindent
 From Figure 6, as the sample size  $n$ increases, the power of the tests for normality and for symmetry increases to one. For $n=50$ the $Z_A$ test is more powerful than the other tests for normality. The tests KS, D (for a sample size of $n\le200$) and BS (for a sample size of $n=50$) have a poor performance. The power of the KS test is less than the nominal level for a sample size of $n\le100$. The new proposed test is more powerful than the CM test for symmetry.
 \bigskip

 Let $X$ be gamma distribution with parameters $a$ and $\theta$. The r.v.\ $X$ is asymmetric for all $a$ and $\theta$ (notice that the test of symmetry of Cabilio and Masaro cannot be used). For large values of $a$, via Central Limit Theorem, the density function of $X$ is close to a normal distribution.

 \noindent
 The parameter $\theta$ does not affect the critical region of the tests (so we fix $\theta=1$), in contrast to the parameter $a$.
 In order to investigate the level of rejection of the tests, as $a$ increases, we choose to simulate the gamma distribution for $\theta=1$ and $a=1$, $10$, $30$, $50$.
 \medskip

 \noindent
 {\footnotesize {\bf Table 7.} Mean and mean square error (MSE) for the three estimators calculated from $10^4$ random samples of size $n$ simulated from the exponential distribution with $\theta=1$ (i.e.\ $\delta=0$, $\beta=1$, $\gamma=1$). The last row gives the empirical size, $p_\delta$, of the test for $\delta=0$ with a nominal significance level of $\alpha=0.05$.}

 \centerline{\footnotesize
 \begin{tabular}{@{\hspace{.5ex}}l
                 @{\hspace{1.3ex}}r
                 @{\hspace{1.3ex}}r
                 @{\hspace{1.3ex}}r
                 @{\hspace{1.3ex}}r
                 @{\hspace{1.3ex}}r
                 @{\hspace{1.3ex}}r
                 @{\hspace{1.3ex}}r
                 @{\hspace{.5ex}}}
   \hline
          $n$                       &
          \multicolumn{1}{c}{$50$}  &
          \multicolumn{1}{c}{$100$} &
          \multicolumn{1}{c}{$150$} &
          \multicolumn{1}{c}{$200$} &
          \multicolumn{1}{c}{$300$} &
          \multicolumn{1}{c}{$400$} &
          \multicolumn{1}{c}{$500$} \\ 
          $\widehat{\delta}_n$      & $-0.2857$ & $-0.1855$ & $-0.1453$ & $-0.1201$ & $-0.0910$ & $-0.0753$ & $-0.0650$  \\
   [-.2ex]$MSE(\widehat{\delta}_n)$ & $ 0.1131$ & $ 0.0521$ & $ 0.0349$ & $ 0.0258$ & $ 0.0172$ & $ 0.0133$ & $ 0.0110$  \\
          $\widehat{\beta}_n$       & $ 1.2820$ & $ 1.1912$ & $ 1.1599$ & $ 1.1371$ & $ 1.1128$ & $ 1.0938$ & $ 1.0820$  \\
   [-.2ex]$MSE(\widehat{\beta}_n)$  & $ 0.5744$ & $ 0.2271$ & $ 0.1433$ & $ 0.1021$ & $ 0.0637$ & $ 0.0450$ & $ 0.0342$  \\
          $\widehat{\gamma}_n$      & $ 1.2755$ & $ 1.1689$ & $ 1.1314$ & $ 1.1105$ & $ 1.0853$ & $ 1.0696$ & $ 1.0613$  \\
   [-.2ex]$MSE(\widehat{\gamma}_n)$ & $ 0.3895$ & $ 0.1510$ & $ 0.0936$ & $ 0.0674$ & $ 0.0419$ & $ 0.0309$ & $ 0.0243$  \\
          $p_{\delta}$              & $ 0.0717$ & $ 0.0385$ & $ 0.0266$ & $ 0.0212$ & $ 0.0170$ & $ 0.0132$ & $ 0.0122$  \\
   \hline
 \end{tabular}}
 \medskip

 \noindent
 From Table 7, the test for $\delta=0$ is liberal for a sample size of $n=50$ and is conservative for $n\ge100$.
 \smallskip

 \noindent
 From Figure 7, as the sample size  $n$ increases, the power of the tests for normality and of the new proposed test for symmetry increases to one. The test KS for a sample size of $n=50$ has a poor performance.

 \centerline{\includegraphics[scale=.99]{Figure.7}}

  \noindent
 {\footnotesize {\bf Figure 7.} Empirical power of the tests for normality (a) and symmetry (b), for a nominal size of 0.05, calculated from $10^4$ random samples of size $n$ simulated from the exponential distribution with $\theta=1$. The tests are the same as those described in the caption to Figure 2.}
 \medskip

 \noindent
 {\footnotesize {\bf Table 8.} Mean and mean square error (MSE) for the three estimators calculated from $10^4$ random samples of size $n$ simulated from the gamma distribution with $a=10, \theta=1$ (i.e.\ $\delta=0$, $\beta=1$, $\gamma=10$). The last row gives the empirical size, $p_\delta$, of the test for $\delta=0$ with a nominal significance level of $\alpha=0.05$.}

 \centerline{\footnotesize
 \begin{tabular}{@{\hspace{.5ex}}l
                 @{\hspace{1.3ex}}r
                 @{\hspace{1.3ex}}r
                 @{\hspace{1.3ex}}r
                 @{\hspace{1.3ex}}r
                 @{\hspace{1.3ex}}r
                 @{\hspace{1.3ex}}r
                 @{\hspace{1.3ex}}r
                 @{\hspace{.5ex}}}
   \hline
          $n$                       &
          \multicolumn{1}{c}{$50$}  &
          \multicolumn{1}{c}{$100$} &
          \multicolumn{1}{c}{$150$} &
          \multicolumn{1}{c}{$200$} &
          \multicolumn{1}{c}{$300$} &
          \multicolumn{1}{c}{$400$} &
          \multicolumn{1}{c}{$500$} \\ 
          $\widehat{\delta}_n$      & $-0.1075$ & $-0.0618$ & $-0.0436$ & $-0.0348$ & $-0.0259$ & $-0.0194$ & $-0.0153$  \\
   [-.2ex]$MSE(\widehat{\delta}_n)$ & $ 0.0251$ & $ 0.0113$ & $ 0.0071$ & $ 0.0055$ & $ 0.0038$ & $ 0.0029$ & $ 0.0023$  \\
          $\widehat{\beta}_n$       & $ 1.0292$ & $ 1.0202$ & $ 1.0140$ & $ 1.0085$ & $ 0.9986$ & $ 1.0065$ & $ 1.0042$  \\
   [-.2ex]$MSE(\widehat{\beta}_n)$  & $ 0.5144$ & $ 0.2239$ & $ 0.1483$ & $ 0.1053$ & $ 0.0681$ & $ 0.0510$ & $ 0.0405$  \\
          $\widehat{\gamma}_n$      & $10.8993$ & $10.5233$ & $10.3646$ & $10.2952$ & $10.1970$ & $10.1635$ & $10.1317$  \\
   [-.2ex]$MSE(\widehat{\gamma}_n)$ & $ 8.6087$ & $ 3.8083$ & $ 2.5283$ & $ 1.8026$ & $ 1.2040$ & $ 0.9187$ & $ 0.7057$  \\
          $p_{\delta}$              & $ 0.0808$ & $ 0.0582$ & $ 0.0478$ & $ 0.0433$ & $ 0.0384$ & $ 0.0322$ & $ 0.0357$  \\
   \hline
 \end{tabular}}
 \medskip

 \centerline{\includegraphics[scale=.99]{Figure.8}}

  \noindent
 {\footnotesize {\bf Figure 8.} Empirical power of the tests for normality (a) and symmetry (b), for a nominal size of 0.05, calculated from $10^4$ random samples of size $n$ simulated from the gamma distribution with $a=10, \theta=1$. The tests are the same as those described in the caption to Figure 2.}
 \medskip

 \noindent
 From Table 8, the test for $\delta=0$ maintains the nominal level well for a sample size of $n=150$, is liberal for $n\le100$ and conservative for $n\ge200$.
 \smallskip

 \noindent
 From Figure 8, as the sample size  $n$ increases, the power of the tests new proposed, BS, AD, CvM, $Z_A$ and $Z_C$ for normality and of the new proposed for symmetry increases to one; also, the power of the tests KS and D for normality increases. The $Z_A$, $Z_C$ and new proposed tests for normality are more powerful than the other tests for normality. The power of KS test for a sample size of $n\le200$ is less than the nominal level. All the tests (new proposed for normality, BS, AD, CvM, $Z_A$ and $Z_C$ for a sample size of $n=50$; D for $n\le200$; KS for $300\le n\le 500$; new proposed for symmetry for $n=50$) have a poor performance.
 \medskip

 \noindent
 {\footnotesize {\bf Table 9.} Mean and mean square error (MSE) for the three estimators calculated from $10^4$ random samples of size $n$ simulated from the gamma distribution with $a=30, \theta=1$ (i.e.\ $\delta=0$, $\beta=1$, $\gamma=30$). The last row gives the empirical size, $p_\delta$, of the test for $\delta=0$ with a nominal significance level of $\alpha=0.05$.}

 \centerline{\footnotesize
 \begin{tabular}{@{\hspace{.5ex}}l
                 @{\hspace{1.3ex}}r
                 @{\hspace{1.3ex}}r
                 @{\hspace{1.3ex}}r
                 @{\hspace{1.3ex}}r
                 @{\hspace{1.3ex}}r
                 @{\hspace{1.3ex}}r
                 @{\hspace{1.3ex}}r
                 @{\hspace{.5ex}}}
   \hline
          $n$                       &
          \multicolumn{1}{c}{$50$}  &
          \multicolumn{1}{c}{$100$} &
          \multicolumn{1}{c}{$150$} &
          \multicolumn{1}{c}{$200$} &
          \multicolumn{1}{c}{$300$} &
          \multicolumn{1}{c}{$400$} &
          \multicolumn{1}{c}{$500$} \\ 
          $\widehat{\delta}_n$      & $-0.0919$ & $-0.0488$ & $-0.0336$ & $-0.0257$ & $-0.0173$ & $-0.0146$ & $-0.0121$  \\
   [-.2ex]$MSE(\widehat{\delta}_n)$ & $ 0.0214$ & $ 0.0088$ & $ 0.0056$ & $ 0.0042$ & $ 0.0028$ & $ 0.0022$ & $ 0.0017$  \\
          $\widehat{\beta}_n$       & $ 1.0068$ & $ 0.9964$ & $ 0.9951$ & $ 0.9996$ & $ 1.0031$ & $ 1.0003$ & $ 0.9946$  \\
   [-.2ex]$MSE(\widehat{\beta}_n)$  & $ 1.1458$ & $ 0.5417$ & $ 0.3446$ & $ 0.2526$ & $ 0.1655$ & $ 0.1241$ & $ 0.1025$  \\
          $\widehat{\gamma}_n$      & $32.0871$ & $31.1106$ & $30.7382$ & $30.6222$ & $30.4259$ & $30.3779$ & $30.3113$  \\
   [-.2ex]$MSE(\widehat{\gamma}_n)$ & $60.6143$ & $28.0028$ & $18.1469$ & $13.6622$ & $ 8.8630$ & $ 7.0383$ & $ 5.3907$  \\
          $p_{\delta}$              & $ 0.0773$ & $ 0.0620$ & $ 0.0477$ & $ 0.0486$ & $ 0.0413$ & $ 0.0415$ & $ 0.0373$  \\
   \hline
 \end{tabular}}
 \medskip

 \centerline{\includegraphics[scale=1]{Figure.9}}

 \noindent
 {\footnotesize {\bf Figure 9.} Empirical power of the tests for normality (a) and symmetry (b), for a nominal size of 0.05, calculated from $10^4$ random samples of size $n$ simulated from the gamma distribution with $a=30, \theta=1$. The tests are the same as those described in the caption to Figure 2.}
 \medskip

 \noindent
 From Table 9, the test for $\delta=0$ maintains the nominal level well for a sample size of $150\le n\le200$, is liberal for $n\le100$ and conservative for $n\ge300$.
 \smallskip

 \noindent
 From Figure 9, as the sample size  $n$ increases, the power of the tests increases. For $n\le150$ the $Z_A$ and $Z_C$ tests are more powerful than the other tests for normality and for $n\ge300$ the new proposed test for normality is more powerful than the other tests for normality. The power of KS test for a sample size of $n\le500$ is less than the nominal level. The tests new proposed for normality, $Z_A$, $Z_C$ (for $n\le100$), BS (for $n\le150$) AD, CvM (for $n\le200$), D (for $n\le500$) and new proposed for symmetry (for $n\le150$) have a poor performance.
 \medskip

 \noindent
 {\footnotesize {\bf Table 10.} Mean and mean square error (MSE) for the three estimators calculated from $10^4$ random samples of size $n$ simulated from the gamma distribution with $a=50, \theta=1$ (i.e.\ $\delta=0$, $\beta=1$, $\gamma=50$). The last row gives the empirical size, $p_\delta$, of the test for $\delta=0$ with a nominal significance level of $\alpha=0.05$.}

 \centerline{\footnotesize
 \begin{tabular}{@{\hspace{.5ex}}l
                 @{\hspace{1.3ex}}r
                 @{\hspace{1.3ex}}r
                 @{\hspace{1.3ex}}r
                 @{\hspace{1.3ex}}r
                 @{\hspace{1.3ex}}r
                 @{\hspace{1.3ex}}r
                 @{\hspace{1.3ex}}r
                 @{\hspace{.5ex}}}
   \hline
          $n$                       &
          \multicolumn{1}{c}{$50$}  &
          \multicolumn{1}{c}{$100$} &
          \multicolumn{1}{c}{$150$} &
          \multicolumn{1}{c}{$200$} &
          \multicolumn{1}{c}{$300$} &
          \multicolumn{1}{c}{$400$} &
          \multicolumn{1}{c}{$500$} \\ 
          $\widehat{\delta}_n$      & $ -0.0873$ & $-0.0451$ & $-0.0320$ & $-0.0238$ & $-0.0169$ & $-0.0123$ & $-0.0098$  \\
   [-.2ex]$MSE(\widehat{\delta}_n)$ & $  0.0199$ & $ 0.0083$ & $ 0.0053$ & $ 0.0039$ & $ 0.0026$ & $ 0.0019$ & $ 0.0015$  \\
          $\widehat{\beta}_n$       & $  1.0241$ & $ 0.9797$ & $ 0.9932$ & $ 0.9919$ & $ 0.9931$ & $ 0.9981$ & $ 0.9994$  \\
   [-.2ex]$MSE(\widehat{\beta}_n)$  & $  1.8321$ & $ 0.8343$ & $ 0.5527$ & $ 0.4124$ & $ 0.2649$ & $ 0.1994$ & $ 0.1610$  \\
          $\widehat{\gamma}_n$      & $ 53.2245$ & $51.7797$ & $51.3292$ & $50.9261$ & $50.6445$ & $50.4826$ & $50.3935$  \\
   [-.2ex]$MSE(\widehat{\gamma}_n)$ & $163.7599$ & $75.7973$ & $51.7595$ & $36.5868$ & $23.9070$ & $18.0106$ & $14.1179$  \\
          $p_{\delta}$              & $  0.0746$ & $ 0.0542$ & $ 0.0492$ & $ 0.0504$ & $ 0.0462$ & $ 0.0432$ & $ 0.0449$  \\
   \hline
 \end{tabular}}
 \medskip

 \centerline{\includegraphics[scale=1]{Figure.10}}

  \noindent
 {\footnotesize{\bf Figure 10.} Empirical power of the tests for normality (a) and symmetry (b), for a nominal size of 0.05, calculated from $10^4$ random samples of size $n$ simulated from the gamma distribution with $a=50, \theta=1$. The tests are the same as those described in the caption to Figure 2.}
 \medskip

 \noindent
 From Table 10, the test for $\delta=0$ maintains the nominal level well for a sample size of $100\le n\le300$, is liberal for $n=50$ and conservative for $n\ge400$.
 \smallskip

 \noindent
 From Figure 10, as the sample size  $n$ increases, the power of the tests increases. For $n\le150$ the $Z_A$ and $Z_C$ tests are more powerful than the other tests for normality and for $n\ge300$ the new proposed test for normality is more powerful than the other tests for normality. The power of KS test for a sample size of $n\le500$ is less than the nominal level. The tests new proposed for normality, BS, AD, CvM, $Z_A$, $Z_C$ (for $n\le200$), D (for $n\le500$) and new proposed for symmetry (for $n\le200$) have a poor performance.
 \bigskip

 \noindent
 {\small\bf Discrete case} \ In order to investigate the results of Subsection \ref{subsec4.2}, we simulated $10^4$ samples of size $n$ from various Poisson, binomial and negative binomial distributions. The Figure 11 shows the empirical size of the proposed test for Poisson for a nominal significance level of $\alpha=0.05$. The Figures 12 and 13 show the empirical power of this test for binomial and negative binomial distributions respectively (again, for a nominal significance level of 0.05).
 \smallskip

 For the Poisson distribution the parameter $\lambda$ affects the critical region of the test. In order to investigate the size of the test, as $\lambda$ increases, we choose to simulate the Poisson distribution for $\lambda=1,10,30$ and $100$.
 \medskip

 \centerline{\includegraphics[scale=1]{Figure.11}}

 \noindent
 {\footnotesize {\bf Figure 11.}
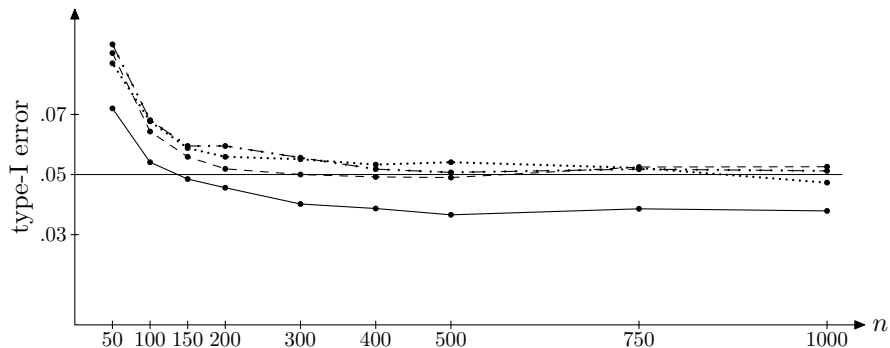
  Empirical type I error calculated from $10^4$ random samples of size $n=50$, $100$, $150$, $200$, $300$, $400$, $500$, $750$, $1000$ simulated from the Poisson distribution with parameter $\lambda=1,10,30,100$ (lines $-$\hspace{-.35ex}$-$, - - -, $\cdots$, - $\cdot$ -) of the test for an underlying Poisson distribution for a nominal significance level of $\alpha=0.05$.}
 \bigskip

 In order to investigate the power of the proposed test for Poisson, we choose to simulate various binomial and negative binomial distributions, see Figures 12 and 13 respectively. Note that all binomial distributions (with parameters $N$, $p$) and all negative binomial distributions (with parameters $r$, $p$, i.e.\ with p.m.f.\ $f(j)={j+r-1\choose j}p^r(1-p)^j$, $j=0,1,\ldots$) satisfy {\rm(\ref{qua2})}.
 \bigskip

 \centerline{\includegraphics[scale=.98]{Figure.12}}

 \noindent
 {\footnotesize {\bf Figure 12.}
  Empirical power calculated from $10^4$ random samples of size $n=50$, $100$, $150$, $200$, $300$, $400$, $500$, $750$, $1000$ simulated from the binomial distribution with parameters $N=10,30$, $100$ (lines $-$\hspace{-.35ex}$-$, - - -, $\cdots$) and $p=0.20,0.35,0.50$, $0.65$ (points \CircleB, \Circle, \SquareB, \Square) of the test for an underlying Poisson distribution for a nominal significance level of $\alpha=0.05$.}
 \medskip

 \centerline{\includegraphics[scale=.98]{Figure.13}}
 \smallskip

 \noindent
 {\footnotesize {\bf Figure 13.}
  Empirical power calculated from $10^4$ random samples of size $n=50$, $100$, $150$, $200$, $300$, $400$, $500$, $750$, $1000$ simulated from the negative binomial distribution with parameters $r$ and $p$, for $r=1$, $10$, $50$, $100$ (lines $-$\hspace{-.35ex}$-$, - - -, $\cdots$, - $\cdot$ -) and $p=0.60$, $0.70$, $0.80$, $0.90$ (points \CircleB, \Circle, \SquareB, \Square) of the test for an underlying Poisson distribution for a nominal significance level of $\alpha=0.05$. The results for the cases $(r=1,p=0.9,n\le500)$, $(r=1,p=0.8,n\le100)$ and $(r=1,p=0.7,n=50)$ are not given, since the estimated matrix $\widehat{\mathbb{D}_{\tau;0}}$ (see Subsection \ref{subsec4.2}) often is singular in these cases.}
 \medskip

 \noindent
 From Figure 11, the size of the proposed test for Poisson is not affected when $\lambda\ge10$. Only in the case $\lambda=1$ this size is different. In the case $\lambda=1$ this test maintains the nominal level well for a sample size of $100\le n\le150$, is liberal for $n=50$ and conservative for $n\ge200$. In the cases $10\le\lambda\le100$ this test maintains the nominal level well for a sample size of $n\ge150$ and is liberal for $n\le100$.
 \medskip

 \noindent
 From Figure 12, the parameter $N$ does not affect the power of the test, in contrast to the parameter $p$. As $p$ increases, the power of the test increases for each $n$. Also, as $n$ increases, the power of the test increases to one for each $p$. For the cases $(p=0.2,n\le200)$ and $(p=0.35,n=50)$ the test has a poor performance.
 \medskip

 \noindent
 From Figure 13, the parameter $r$ does not affect the power of the test, except in cases $p=0.8$ and $n\le300$. Only the parameter $p$ affects the power of the test. As $p$ increases, the power of the test decreases for each $n$. Also, as $n$ increases, the power of the test increases for $p=0.9$ and increases to one for $p\le0.8$. For the cases $(p=0.8,n\le150)$ and $(p=0.9,n=750)$ the test has a poor performance.

 \noindent
 Notice that we home also simulated negative binomial distributions for $p\le0.5$, $r=1,10,50,100$ of size $n=50$, $100$, $150$, $200$, $300$, $400$, $500$, $750$, $1000$ and in all cases the observed empirical power is one.
 \medskip

 Similar to the continuous case, we can construct tables for the average of the estimators and for their MSEs. From these tables one can see that, as $n$ increases, the averages of the estimators tend to their true values and the mean square errors decrease to zero.
 \bigskip

 \noindent{\sc Acknowledgments:}
 I am most grateful to Professor Balakrishnan for his careful reading of the paper and for his helpful comments and observations.
   \end{document}